\documentclass[twocolumn,aps]{revtex4-2}

\usepackage{graphicx} 
\usepackage{fancyvrb}
\begin{document}
\title{Pearl necklace knots with fewer vertices than an equivalent FCC lattice knot}
\author{Alexander R. Klotz }
\affiliation{Department of Physics and Astronomy, California State University, Long Beach}
\begin{abstract}
    The pearl necklace number, $N_P$, of a knot is the smallest number of unit spheres required to construct a knot if each sphere is tangent to two neighbors and no sphere overlaps with another. It has been speculated that the pearl necklace number is equal to the minimum number of lattice sites required to embed a knot on a face centered cubic (FCC) lattice $N_L$, implying that a trefoil knot cannot be constructed from fewer than 15 spheres. A large language model (LLM) was prompted to find configurations of knot for which $N_P<N_L$, and this manuscript describes its findings, attempts at validating them, and their implications. No 14-vertex example was found for the trefoil knot, but examples with $N_P=N_L-1$ were found for all knots from 5 to 7 crossings as well as $8_{19}$, and $10_{124}$. One example, $8_1$, could be constructed with $N_P=N_L-2$. 
\end{abstract}

\maketitle

\section{Introduction}

A ``pearl necklace knot'' is an equilateral polygonal knot constructed from spheres of equal radius (taken to be 1) such that each sphere is tangent to its two neighbors, with no overlaps between non-adjacent spheres. The pearl necklace number, $N_P$, is the minimum number of spheres required to construct a given knot. Exact values for $N_P$ are unknown. A lower bound on $N_P$ is the stick number of the knot and the upper bound is the minimum number of sites required to embed the knot on a face-centered cubic (FCC) lattice. Values of the FCC lattice number, $N_L$, have been computed for up to 10 crossings, and the trefoil knot has $N_L=15$ \cite{van2011minimal}. Prior to this work, it was unknown whether knots could be constructed such that $N_P<N_L$. In a recent survey of open problems in low-dimensional topology \cite{baykur2026k3}, a question was posed as to whether a trefoil knot can be constructed from 14 or fewer tangent spheres. 

Most of the work on the topic has been carried out by Maehara and collaborators \cite{maehara1999knotted, maehara2000knotted, maehara2007configurations}. They have established a lower bound for the pearl necklace number of the trefoil knot (at 11) by constraining it to the minimum distance (MD) energy of a knot, which increases as the line segments of polygonal knots approach each other. They have also considered the pearl necklace number for knots that lie on a table or in a display case, with constraints on the z-coordinates of the vertices. Recent work has also examined knots constructed from spheres without equal radius \cite{alfonsin2024links}.

Recently, success has been established in prompting large language models (LLMs) to prove or find counterexamples to various open problems in mathematics. Intrigued by these reports, I began to prompt an LLM known as ChatGPT (version GPT-5.6 Pro) about various questions in knot theory. When addressing the question of the pearl knot number, it began to find several results for knots with $N_P<N_L$ and I began to prompt it further. In this manuscript, I will describe the results found by ChatGPT, my efforts to validate and refine them, and their implications. A rendering of the transcript may be found in the ancillary files.
\begin{figure}
    \centering
    \includegraphics[width=1\linewidth]{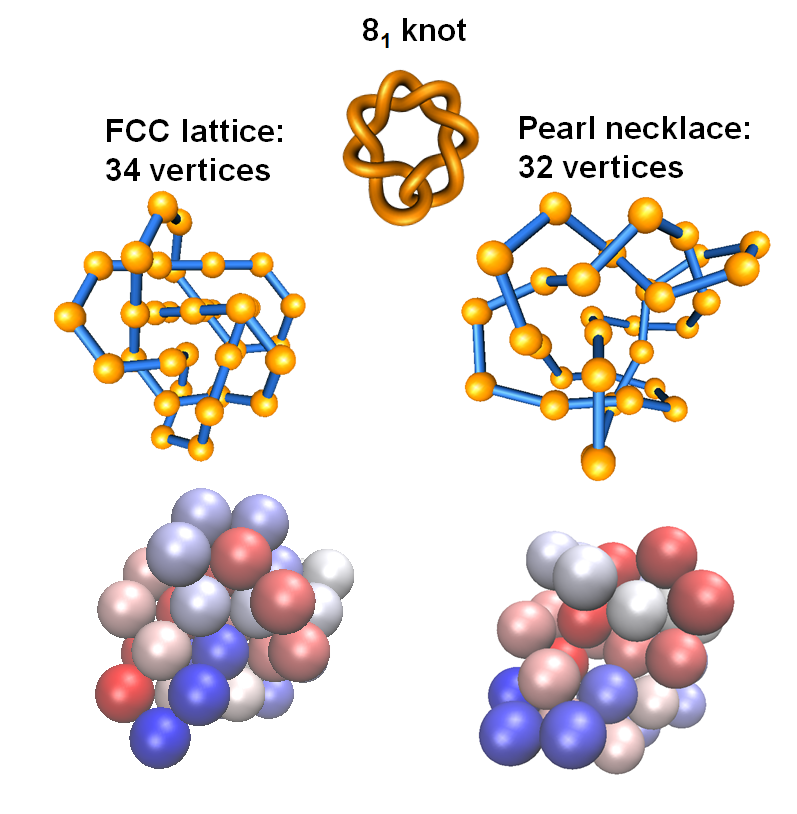}
    \caption{An $8_1$ knot embedded on an FCC lattice with 34 sites (left) and constructed from 32 tangent spheres (right). Top images show a ball-and-stick visualization, bottom images show space-filling spheres color-coded by index.} 
    \label{fig:f2}
\end{figure}

\section{Results and Validation}
I prompted ChatGPT to find a 14-sphere trefoil configuration. This was unsuccessful, with the minimum distance between non-adjacent spheres being 1.98579138665. I then prompted it to find a 32-vertex configuration for the $10_{124}$ knot (a non-alternating torus knot), for which the FCC lattice number is 33. It was successful, and I repeated the prompt for the $8_{19}$ knot (also a non-alternating torus knot with $N_L=28$) where it was again successful. I initially suggested non-alternating torus knots as I was aware from my work with the ropelength problem that these knots typically have the smallest ropelength at a given crossing number, and thus might require the fewest amount of spheres as pearl necklaces. I prompted it for the $4_1$ knot (unsuccessful), and in each of these cases I was supplying it with the number of vertices to target. After the $4_1$, I prompted it to look at Andrew Rechnitzer's website \cite{ubcMinimalKnots}, which contains the lattice coordinates for all knots up to 10 crossings, and then find solutions for all knots with up to 7 crossings with one fewer vertex, which were all successful. The methods described by the LLM involve removing a vertex from the lattice configuration and evolving the resulting configuration towards an equilateral state while preserving topology. The results are summarized in Table I, and coordinates are provided in Appendix A.

Afterwards, I prompted the LLM to check knots with 8-10 crossings for examples where more than one vertex could be removed. It found the $8_1$ (a twist knot, Figure 1) as an example where two vertices could be removed, from 34 down to 32, but no others. This is likely the most significant finding. This search was non-exhaustive (it checked $8_{19}$, $9_1$, and 10-crossing knots 1, 20, 55, 82, 100, 124, 127, 135, 146, 153),  but other examples may exist with 8-10 crossings.

\begin{table}[ht]
\centering
\caption{Knots found with $N_P<N_L$.}
\label{tab:knot_clearance}
\begin{tabular}{c c c c}
\hline
Knot
& $N_L$
& Best $N_P$
& Minimum clearance \\
\hline
$5_1$      & 22 & 21 & $0.000741768$ \\
$5_2$      & 23 & 22 & $0.029976777$ \\
$6_1$      & 27 & 26 & $0.029980437$ \\
$6_2$      & 27 & 26 & $0.029920813$ \\
$6_3$      & 28 & 27 & $0.029972539$ \\
$7_1$      & 29 & 28 & $0.029959263$ \\
$7_2$      & 30 & 29 & $0.029900448$ \\
$7_3$      & 30 & 29 & $0.029860976$ \\
$7_4$      & 30 & 29 & $0.009917403$ \\
$7_5$      & 31 & 30 & $0.029966925$ \\
$7_6$      & 31 & 30 & $0.029969891$ \\
$7_7$      & 31 & 30 & $0.029952571$ \\
\textbf{$8_1$ }     & \textbf{34} & \textbf{32} & \textbf{$0.014334921$} \\
$8_{19}$   & 28 & 27 & $0.020367894$ \\
$10_{124}$ & 33 & 32 & $0.012391365$ \\
\hline
\end{tabular}
\end{table}

Three constraints must be satisfied for these configurations to be valid: the knot topology must not change, adjacent vertices must be separated by a distance of exactly 2, and non-adjacent vertices must be separated by at least 2. An evaluation of the Alexander polynomial of each configuration was used to verify that it had not undergone a topological change during perturbation. It is straightforward to verify that a configuration is non-overlapping as long as no pair of non-adjacent vertices are within distance 2 of each other, provided the minimum pairwise distance is not sufficiently close to 2 to suspect rounding errors. Many of the minimum clearances (the smallest distance above 2 between non-adjacent spheres) are close to 0.0299. This likely arises because target was only a non-overlapping configuration, not one that maximizes this parameter, but may hint at some deeper geometric constraints. 

Verifying that the configurations are equilateral is more challenging. The initial configurations provided by the LLM typically had deviations from edge length of two on the order of $10^{-13}$. While small, that is several orders of magnitude larger than the numerical precision of a float variable. The configurations were refined according to a human-created procedure as follows:
\begin{itemize}
    \item For each vertex, compute the $i$th unit vector $\hat{n}_i$ between pairs of vertices at $\vec{r}_i$  and $\vec{r}_{i+1}$, defined cyclically.
    \item Displace each subsequent point for $i>1$ by $2\hat{n}_i$ from the previous such that $\vec{r}_{i+1}=\vec{r}_i+2\hat{n}_i$. This has the effect of collating all the equilateral deviation between the final and first vertices. 
    \item Define point $\vec{P}$ as the midpoint between $\vec{r}_1$ and $\vec{r}_{N_P-1}$, with $|v|$ being the distance between them and unit vector $\hat{v}$ along the direction between them.
    \item Find the displacement between the last vertex $\vec{r}_{N_P}$ and the plane passing through $\vec{p}$ with normal $\hat{v}$ as $\delta=\hat{v}\cdot(\vec{r}_{N_P}-\vec{P})$ and subtract that distance from the last vertex along $\hat{v}$. The final vertex now lies in the plane directly between its neighbors.
    \item Define the unit vector $\hat{v}_2$ between the new location of $\vec{r}_{N_P}$ and $\vec{P}$. 
    \item Find the distance $h$ between the final vertex from $\vec{P}$ that would allow its two adjacent edges to be 2, $h=\sqrt{4-|v|^2/4}$ .
    \item Place the final vertex at position $\vec{r}_{N_P}=\vec{b}+h\hat{v}_2$.
    
\end{itemize}
This refinement ignores the rotational symmetry for the location of the final vertex in the bisecting plane, but that has not been an issue. After such refinement, the deviation from equilateral edges is effectively zero or within the uncertainty of a float variable, $2^{-52}$. For example, the configuration shown for the $8_1$ knot has seven edges with length $2-2.22\cdot10^{-16}$ and three with $2+4.44\cdot10^{-16}$, where $2.22\cdot10^{-16}\approx2^{-52}$, and the rest exactly 2. A final check is necessary to determine whether the refinement has left the knot in a non-overlapping configuration, but the displacements imposed during this procedure are a billionth of the typical clearance between non-adjacent spheres.  

To further verify whether the configurations were truly acceptable, the LLM was prompted to generated algebraic expressions for the vertices such that the equilateral and non-overlapping could be verified exactly. It output a scheme, described in Appendix B, by which each edge vector could be described as a harmonic function of a rational multiple of $\pi$, with the fractions having numerators and denominators of order 1000. The knot could be constructed as the cumulative sum of those harmonic vectors, with a procedure for finding the positions of the final two coordinates to ensure equilateral closure between the first and last vertex. A table of fractions for the $5_1$ knot is given in the Appendix.


Sufficiently complex knots may allow configurations in which the vast majority of vertices lie on a lattice but those near the surface of the knot may be perturbed and removed. To determine whether this is feasible, the LLM was prompted to find pearl necklace configurations in which as many vertices as possible occupied lattice sites. It found a configuration for the $6_1$ knot in which 16 of 26 sites lay on an FCC lattice (Figure 2), an $8_1$ in which 26 of 33 sites lay on the lattice, and a $10_{124}$ in which 22 of 32 sites lay on the lattice. For twist knots, it was apparent that the twisted regions of the knot could remain on the lattice while the clasp was perturbed from it. Coordinates for these configurations are found at the end of Appendix A.

\begin{figure}
    \centering
    \includegraphics[width=0.8\linewidth]{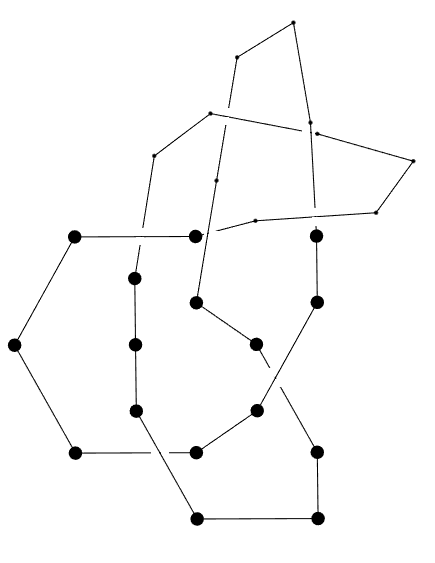}
    \caption{A $6_1$ pearl necklace knot in which 16 of 26 sites (with larger markers) lie on an FCC lattice. 27 sites are required to fully embed the same knot on the lattice. } 
    \label{fig:f2}
\end{figure}

Most of the methods used by the LLM appear to be based on gradient descent towards an equilateral configuration, supported by impressive computational power. The LLM was prompted to produce scripts that would replicate its results. The first script, in Python, was too slow to validate on my computer. The second, in MATLAB, was tested on a few knots but did not produce non-overlapping configurations during my attempts at validation. These scripts are provided in the ancillary material. Another strategy for generating such knots would be to perturb an FCC knot using length- and topology-preserving crankshaft moves until three subsequent vertices form an equilateral triangle, after which the central vertex can be removed. The LLM was unable to find a sequence of moves that would allow this.

\section{Discussion}

These examples show that there are knots for which $N_P<N_L$, although the 14-pearl trefoil conjecture remains open. There are two potential reasons why these results would \textit{not} suggest that the pearl necklace number can be smaller than the FCC lattice number. The first is that the numerical imprecision is actually significant and these knots aren't truly equilateral. I deem this unlikely after the knots survive the refinement procedure, and exact lengths can be proscribed with the harmonic parameterization. The second is that the FCC lattice number is actually smaller than that reported by Janse van Rensburg and Rechnitzer \cite{van2011minimal}. There may be some examples at higher crossing number where this is the case, but I do not suspect it to be so for the simpler knots (those with up to 8 crossings, appearing in their paper). Their methods are robust and their computations were thorough, stating that billions of trial BCACF moves were attempted. A difficult alternative to BCACF sampling would be to check exhaustively every distinct self-avoiding polygon of length 21 (for the $5_1$ knot with $N_L=22$) to determine if any are knotted. 

We may ask whether sufficiently complex knots still have $N_P$ only deviate by 1 or 2 from $N_L$, or if the gulf increases with crossing number. Extending this work to larger crossing number is limited by the fact that $N_L$ has ony been computed up to 10 crossings. The knots that have been minimized on lattices thus far are simple enough that most occupied sites do not have all 12 neighboring sites occupied. For very large knots (as the crossing number tends to infinity), we would expect the vast majority of lattice sites to be in the ``bulk'' of the knot. In such cases, the only vertices that could be perturbed (initially) would be on the surface. If we assume that the minimum lattice number of a knot depends similarly on crossing number as ropelength does (which is reasonable given that lattice embeddings are used to derive ropelength bounds \cite{diao}), we expect asymptotically that $N_L\propto C$ for alternating knots and $N_L\propto C^{3/4}$ for non-alternating knots. The number of sites on the surface of the knot we expect to grow as $N_L^{2/3}$. If many surface sites can be perturbed until they are able to be removed, we can hypothesize in the absence of data that:
\begin{equation}
    (N_L-N_p)\propto C^{2/3} \ \ \ \ \  \ \ \text{(alternating)}
\end{equation}

\begin{equation}
    (N_L-N_p)\propto C^{1/2} \ \ \ \ \  \ \ \text{(non-alternating)}
\end{equation}

Existing lower bounds for $N_P$ are based on constraining it to the MD energy. The pearl necklace knots did not minimize MD energy for their vertex count, for example the $5_1$ pearl necklace knot had an MD energy of near 144, which could be minimized to 133. A potentially interesting observation is that the number of vertices that minimizes MD energy for a knot appears to be close to $N_L$ and $N_P$. A pathway towards developing a proof that certain numbers are minimal may involve deriving the MD energy minimum for equilateral knots and comparing it to the global energy minimum at a given vertex count.

Throughout the rise of LLMs over the past few years, I have been a proponent of human thought and creativity over the convenience of a text generator. However, these systems have reached the point where they may be prompted to make discoveries that elude humans. After arriving at these results by prompting an LLM, I began to wonder whether I could have arrived at these results myself if I had put my mind to it. Prior to this investigation, I had attempted to find such configurations by perturbing lattice knots by crankshaft moves until one vertex could be removed, but this was unsuccessful. After observing the LLM results, I wrote several scripts that would remove a vertex from an FCC knot and attempt to minimize various knot energies (e.g. with penalties for non-equilateral bonds, bond crossing, and sphere overlaps) but was ultimately unsuccessful. Could I have eventually found such configurations without using an LLM? I do not know.

To avoid confirmation bias in the assessment of LLMs as a research tool, it is important to report non-successes as well. Prior to investigating the pearl necklace problem, I had prompted the LLM to construct a non-trivial knot with a trivial Jones polynomial. It is conjectured that no such knot exists. The LLM did not find any examples that would disprove that conjecture. I had also prompted it to find stronger lower bounds for the ropelength of an ideal knot as a function of crossing number. It purported a stronger constant for Buck's three-quarter coefficient of 1.105 \cite{buck1999thickness}, allegedly proving that the ropelength of a knot with crossing number $C$ must be at least $(4\pi/4.724)^{3/4}C^{3/4}\approx2.08C^{3/4}$. I have not evaluated the rigor of this claim, but I have included the output from the LLM for interested readers.

Finally, the LLM that produced these configurations was asked to describe what it had found:

\begin{quote}
\textbf{Write a one paragraph summary of your methods and findings from this discussion as well as their significance. Assume the reader has an interest and some background in physical knot theory.
}  \\  We investigated discrete physical-knot realizations built from congruent radius-1 spheres whose centers form a closed equilateral polygon with consecutive center spacing exactly 2 and all non-neighbor separations at least 2, using FCC lattice embeddings and standard torus-knot parametrizations as topology-controlled seeds and then relaxing the vertices off-lattice by constrained nonlinear optimization. Candidate configurations were evaluated by edge-length error, minimum non-neighbor sphere clearance, nonadjacent-edge separation, and preservation of knot type under explicit deformations, with successful examples further converted into exact algebraic forms by replacing numerical edge directions with rational-angle trigonometric vectors and solving the final closure constraints exactly. This approach produced off-lattice constructions that improve on several FCC lattice vertex counts, including one-below-FCC examples for all prime knots through seven crossings that we tested, a 27-sphere realization of $8_{19}$, a 32-sphere realization of $10_{124}$, and a two-below-FCC construction for $8_1$; by contrast, repeated 14-vertex trefoil searches from very different initial conditions converged to the same apparent optimum with minimum non-neighbor distance about ($1.98579138665<2$), suggesting a genuine geometric obstruction or highly structured local/global optimum. The significance is that lattice stick/sphere counts can substantially overestimate the true unconstrained physical complexity of a knot, while the repeated near-optimal failures—especially for the 14-vertex trefoil—point toward contact-graph and algebraic methods that may eventually turn these numerical observations into exact bounds or proofs.

\end{quote}

\section{Acknowledgements}
Access to ChatGPT Pro was provided by the California State University system, an expenditure opposed by the author despite its utility in this manuscript. Discussions with Agnese Barbensi were vital in convincing me that I was not suffering AI-induced hallucinations. In addition to finding the actual results of this research, the LLM was used to produce LaTeX tables but was not involved in the writing of the manuscript.

\section{Ancillary Files}
The ancillary files include an attempted HTML rendering of the conversation with ChatGPT that this paper discusses, two scripts generated by ChatGPT for replicating its results, and a screenshot from a shorter conversation with ChatGPT about ropelength bounds.

\bibliographystyle{unsrt}
\bibliography{pearl}

\clearpage

\appendix
\section{Numerical Coordinates for Off-Lattice Sphere-Knot Configurations}
\label{app:coordinates}

The following coordinate lists give the centers of unit-radius spheres, determined by the LLM. The resulting knots are equilaterall to approximately 12 decimal places; if 16-digit precision is desired then the refinement procedure in the text may be followed.

\tiny
\textbf{$5_1$: 21 vertices}

\begin{verbatim}
 3.618507289129,  0.030051788303,  0.433449301574
 2.087353466509,  1.189770078362,  0.990783545256
 0.366673509445,  0.948520480063,  0.000299261343
-0.744095626595,  2.284975956868, -0.989689973416
-2.657294038455,  2.456428837867, -0.432665964550
-3.167198119193,  0.887511982431,  0.698040774253
-1.594753065537, -0.332356751805,  0.896375905957
-1.037114996792, -1.160278433623, -0.836709629458
-0.168002189633, -2.945431143957, -0.596142114573
 1.553065121170, -2.722156363754,  0.397877847102
 1.834068074356, -0.787778064477,  0.821219240764
 2.085787084637,  0.331858170856, -0.816783274617
 1.708857463053,  2.287911025834, -0.638510200641
 0.274289206620,  2.842086648745,  0.640123170035
-1.320383450218,  1.648106215109,  0.817407804529
-1.886884326119,  0.651186760675, -0.821258836059
-2.980178920339, -0.969440306424, -0.399032783293
-1.857033290712, -2.292857354192,  0.594501922858
-0.013202229765, -1.556628408939,  0.835860491920
 0.956691924369, -1.318134303479, -0.896890468673
 2.940744667913, -1.473264293569, -0.698221404623
\end{verbatim}

\textbf{$5_2$: 22 vertices}

\begin{verbatim}
 0.034903183398,  0.990244692123,  1.520449530016
 1.469083484935,  1.243140912984,  2.891277879650
 1.355927731432,  2.627598094849,  4.330192083843
 2.782177950451,  4.023084505702,  4.465941816341
 3.567290417535,  4.702919383970,  2.756724758060
 2.699381303575,  4.015416253522,  1.091169830857
 1.628249801583,  5.017747552137, -0.268245217337
 0.066999715497,  5.992541663459,  0.514234564657
-0.108097319785,  5.615477579356,  2.470548369933
 1.327373444078,  4.369835023156,  3.093282910711
 2.643065692686,  2.895726344576,  2.783512817438
 3.033725322071,  1.497695963964,  4.159338427337
 1.423324542185,  0.482744032358,  4.772919640115
-0.207812472987,  1.368813166602,  4.028424155189
 0.107227799414,  2.717082069192,  2.585191466718
 1.499398504867,  2.357791955549,  1.194949998082
 2.780549950165,  2.565012642185, -0.326795894386
 3.955340547396,  4.174779475161, -0.495668023850
 3.345531270405,  5.829240014022,  0.448203244411
 1.738358304219,  5.719273746641,  1.633489065130
 0.558629374186,  4.156351616863,  1.226654863158
-0.344489043547,  2.586932130884,  0.377361274003
\end{verbatim}

\textbf{$6_1$: 26 vertices}

\begin{verbatim}
 1.899370116328, -0.246751688529,  3.280951416009
 2.848088774639,  1.413853394517,  2.695867092571
 4.510355147233,  2.459040904297,  2.315796604444
 4.307415578539,  4.237283704120,  1.423233958835
 2.990663988166,  4.367979628288, -0.076460618290
 1.536830655927,  5.740860483440, -0.116043438120
 0.133899732695,  5.829069662454,  1.306631495018
 0.040725822452,  4.408089697954,  2.710958580472
 1.355129671745,  2.900718390142,  2.697757771777
 1.308914496687,  1.403671731005,  1.372346342084
 2.608825243960, -0.116232513423,  1.383443048441
 3.892915267452, -0.314124027953,  2.903953410144
 3.229875264901,  0.638539808015,  4.532698193914
 1.469512474555,  1.520711357583,  4.182132645838
-0.048026847330,  1.445311227157,  2.881598437575
-0.192210621358,  2.770125461557,  1.390264040587
 1.236273883542,  4.137942096070,  1.092758286678
 2.628949703577,  5.489380728812,  1.576563704408
 4.154200479802,  5.968975089365,  0.375064924155
 5.354126451456,  4.515392805141, -0.293722191109
 4.312219088967,  2.827568977880, -0.037441532221
 2.811608901658,  2.871989673394,  1.283995514162
 3.151022057655,  3.786960026218,  3.029741189508
 2.945639719447,  2.792137819300,  4.752572252548
 2.209685220937,  1.481578103580,  6.071967647285
 1.527848169009, -0.207343538003,  5.245746173822
\end{verbatim}

\textbf{$6_2$: 26 vertices}

\begin{verbatim}
 1.095751160782, -0.260375493239,  1.825051523589
 1.433329794707,  0.517940524656,  3.636201667063
 0.891540233423,  1.878892150722,  4.997919237925
 1.062761155103,  3.450098490426,  6.223477899944
 1.546457525292,  5.333312436022,  5.754923680481
 2.917353649998,  5.437049958786,  4.302380276091
 2.715400576810,  4.493929224598,  2.550313428914
 2.455163208135,  3.410619322207,  0.889375089503
 1.750314149298,  1.575552528091,  1.257773875318
 0.109261644644,  1.421230531697,  2.390527775510
-0.525697806732,  0.688967014725,  4.139989588988
 0.722732891921, -0.080859022467,  5.499691270653
 2.597071747192,  0.612491377816,  5.577915806749
 2.821872841110,  2.304927887353,  4.536223371642
 4.026935619568,  3.779527443844,  3.925175216443
 4.440291785965,  5.482302176832,  2.960968346315
 2.772225339696,  6.522652238293,  2.593237357075
 1.160754604954,  5.567131458274,  3.293338856875
 1.588786132309,  3.902348488971,  4.315731480434
 3.040550249556,  3.960626423055,  5.690134852910
 4.552686613505,  2.660816240234,  5.535418643486
 4.618637430026,  1.849832134015,  3.708411454321
 3.101363735732,  2.401510763473,  2.527932218417
 1.231721358819,  3.110203610005,  2.481115215889
 0.138904188683,  2.739951517701,  0.847510606200
 0.059801535781,  0.778332920291,  0.465679166167
\end{verbatim}

\textbf{$6_3$: 27 vertices}

\begin{verbatim}
 0.173446626397,  4.554221121851,  4.378628180113
 1.220000565438,  2.920020861150,  4.862481836299
 3.021545230748,  2.673710535120,  5.695409229404
 4.212403286684,  1.180559164549,  5.101813343993
 4.182585510811,  0.885073780914,  3.123986433392
 3.002317103755,  1.130780044014,  1.528181019320
 3.198727365294,  2.711087079714,  0.318200572337
 4.429130838133,  4.287796944369,  0.327859595326
 3.962385542502,  5.937843490179,  1.357177440148
 2.564390404940,  5.721737355010,  2.771004775326
 1.169710017008,  4.295310472800,  2.628973211135
 1.067873047546,  2.925094401030,  1.175653855895
 1.573062628535,  1.367999180301,  0.026637772090
 3.082761733996,  0.076797619617, -0.204894386921
 4.615452171659,  1.257292049303,  0.302346112033
 4.477141110124,  2.481043404817,  1.878195936773
 2.614346702975,  2.875973327829,  2.489773693567
 2.620034857470,  4.138781918693,  4.040670888012
 1.989343328796,  4.481078907066,  5.907502687455
 1.384960258319,  2.843847232295,  6.884330997470
 1.653723084248,  1.178148094531,  5.810453478452
 2.521115479134,  1.587146085965,  4.055360581680
 4.114966200031,  2.778297785356,  3.853376650304
 4.029340361930,  4.327921209894,  2.591895196627
 2.600848096418,  4.440753433959,  1.196660902114
 1.167552866954,  5.831632282397,  1.302116173033
 0.062048176451,  5.963226914508,  2.963604130279
\end{verbatim}

\textbf{$7_1$: 28 vertices}

\begin{verbatim}
 1.132329610231,  6.528584079107,  6.211000416005
 1.925623986150,  7.413711872669,  4.602512672129
 2.566752454764,  6.233746635690,  3.120407302645
 1.580207683598,  4.542944498659,  2.710629865615
 1.406684453659,  2.981176878359,  1.473386320934
 2.777096105644,  1.524513984175,  1.463133671897
 4.132338249308,  0.808522980301,  0.178347458394
 5.457709547794,  2.204590493352, -0.364226426118
 5.772069557478,  3.822716465615,  0.768401545547
 4.630860806141,  4.024882301442,  2.398362146512
 4.396466302176,  3.204013008113,  4.207017121840
 4.190796875030,  4.370992643306,  5.818183971606
 2.898195712228,  5.856616659201,  5.468750834484
 1.209761053427,  5.520226456401,  4.450895228604
 0.577081073040,  6.307712382635,  2.724748665186
 1.799867709548,  5.881654834765,  1.200521711127
 2.968591449866,  4.259652794866,  1.257038977353
 2.919029207941,  3.021468261065,  2.826893315956
 1.469118438585,  1.656098086435,  3.009986317297
 0.848659637267,  0.925787732813,  1.254515536208
 1.735028176863,  1.929054216604, -0.231353628462
 3.541716639742,  2.745225876029,  0.032730866448
 4.740378570479,  2.191030503639,  1.534755188829
 6.017153100759,  2.668716900578,  2.998197194857
 5.952278632280,  4.465021062596,  3.875170775335
 4.093160454329,  5.196713581850,  3.966306156382
 2.592350185648,  4.042488323102,  4.610769180525
 1.520336533989,  4.566588247129,  6.215794414062
\end{verbatim}

\textbf{$7_2$: 29 vertices}

\begin{verbatim}
 0.579162456958,  4.631373353293,  2.714294056232
 1.756754424950,  4.135373400969,  1.175703401049
 3.231299851243,  2.918470983343,  1.762954873395
 4.635149951207,  2.583556668320,  3.147526033872
 6.042534802461,  3.610881027956,  4.129296139628
 5.549883331474,  5.543682798051,  4.276169709927
 3.788665076065,  5.904497793757,  3.399857667509
 2.567987378027,  4.365894895354,  3.022171185701
 1.408851034300,  2.779991491148,  2.646254533796
 1.121695561138,  1.308706130125,  3.970211149654
 2.809364250272,  0.386646887167,  4.519375194537
 4.396687985426,  0.342393019256,  3.303461976009
 4.354014330387,  1.244891528110,  1.519176292414
 4.128717047717,  2.303693882573, -0.162542681912
 5.421787990384,  3.771918679107, -0.577613138667
 6.098672862173,  4.854831023578,  0.961583965195
 5.028090611220,  4.502248738916,  2.613715732484
 4.092324934295,  4.132393979681,  4.342168813988
 2.517657810949,  5.258380902592,  4.844739107664
 1.789766665359,  6.226470074583,  3.253205936396
 2.841678793699,  5.795878039726,  1.607584055839
 4.110732303053,  4.475945198139,  0.803046620974
 5.359912700061,  2.993965378232,  1.296285051266
 5.905577863312,  1.399455293081,  0.219350922376
 4.334996092979,  0.254769594374, -0.252843480250
 2.747028015009,  0.139051549305,  0.957514614446
 2.666752044264,  1.186788096540,  2.659220786529
 2.743947331206,  2.515307205645,  4.152230453081
 1.207896556408,  3.770637807832,  4.406578829244
\end{verbatim}

\textbf{$7_3$: 29 vertices}

\begin{verbatim}
 1.096526042691,  3.296079831801,  2.271929908764
 2.170834006738,  4.698799931978,  3.209073824576
 2.944736523206,  4.322759874283,  5.014529105111
 2.456383130158,  2.661622369987,  6.015595114583
 2.845723613996,  0.860782977481,  6.793664179669
 4.421018156251, -0.348179321078,  6.555216563186
 5.964682886813,  0.195706394854,  5.405743845970
 5.497883316078,  1.403344014712,  3.881371995054
 5.070064647881,  2.768413897532,  2.483675532241
 6.290008781616,  2.903535304559,  0.904600323161
 7.489166813423,  1.307396987771,  0.784756723280
 7.069234627587, -0.191612555122,  2.040394947298
 5.260105584337, -0.345750961796,  2.879020880867
 3.909356820602,  1.111569935210,  2.651658657927
 3.001126189403,  2.883443946813,  2.840284232960
 1.572318430681,  3.049439273602,  4.229871920046
 1.092348474836,  4.153120912639,  5.827218369987
 2.673755399778,  4.409638989703,  7.024444586658
 4.267764214900,  3.238855061531,  6.727121713599
 4.453977445727,  1.510095931227,  5.738832028306
 4.142543766159, -0.068425247767,  4.550858985390
 3.282452886370, -0.802937573174,  2.901393822881
 4.150372923527, -0.266876228723,  1.181116797491
 5.594289145109,  1.084757444431,  1.478088780424
 6.953225188059,  2.020906701589,  2.608096479289
 6.291682644039,  3.265964764899,  4.026613973520
 4.338451522743,  2.998365492025,  4.363183158027
 2.946181720795,  1.562650114505,  4.380680395975
 1.547724007506,  1.470860065804,  2.953834009323
\end{verbatim}

\textbf{$7_4$: 29 vertices}

\begin{verbatim}
-0.502462877555,  0.685196602755,  1.640444927970
 0.300222571272,  2.129014507228,  2.767869395047
 1.856714506341,  3.127131662504,  2.005573520461
 2.681766454341,  4.639539426585,  3.021404233844
 3.183514648818,  6.026816658277,  4.371852541543
 4.217272060428,  5.321008607733,  5.931718115793
 5.782652689930,  4.081151718993,  6.042796377502
 5.548862681143,  2.247251059282,  6.805784321118
 3.856997074434,  1.431123630908,  6.119106238251
 2.968476550863,  2.410071538053,  4.618375658071
 1.378323555777,  3.413402126897,  3.936656374413
 0.178724710985,  4.117278605396,  2.499464489316
 0.184421993053,  2.897729310023,  0.914325045251
 1.366479492443,  1.323353803094,  1.266577555983
 1.017821708297,  0.126902257307,  2.830847816564
 2.595390885280, -0.219909852930,  4.010253112477
 4.237076600984,  0.903493791124,  4.217212681904
 4.941467404805,  2.634584662085,  4.929362941309
 3.940084973151,  3.406902438154,  6.478801058286
 4.856201148723,  4.675207946089,  7.724647072606
 5.961520975515,  5.977540603693,  6.684359742373
 5.801565005420,  5.600283170458,  4.726787031923
 4.320332767394,  4.383212563823,  4.156971100680
 2.636619347597,  4.294820762232,  5.232748033602
 1.438579195701,  2.819585284528,  5.855949509666
 1.216755765024,  1.444585881765,  4.420619928001
 2.522268698892,  1.457653936963,  2.905540654353
 2.667210567055, -0.142189637982,  1.714116132799
 0.871263015229, -0.600852732208,  0.962983283640
\end{verbatim}

\textbf{$7_5$: 30 vertices}

\begin{verbatim}
 1.446128020728, -0.698615617144,  1.074663208035
 2.387307121994,  0.414056550913, -0.295060320919
 3.932976928993,  1.553771099443,  0.263469531103
 4.911262292926,  1.226779177593,  1.976955460359
 6.392656040307,  1.639957766746,  3.255531237811
 7.159680261035,  2.764609973534,  4.720738137393
 5.958460149612,  4.187169580603,  5.451076694672
 4.184158115208,  4.330201619290,  4.539268531394
 2.664021987057,  3.074361108254,  4.204528688037
 1.356149325789,  2.745080992626,  2.727692980398
 1.664058829422,  1.295726421619,  1.384352716327
 3.261494826783,  0.098863564224,  1.509722309627
 4.882068855889, -0.236453356550,  0.386645117788
 5.959911720168,  1.442892029517,  0.252280560642
 5.624665143698,  3.040192963881,  1.408238802987
 4.263651843673,  3.182421234488,  2.866804508122
 2.713355375786,  4.434735957283,  2.698612068565
 2.844433298834,  5.980257943400,  1.436003511442
 4.733573518746,  6.361249215158,  0.901218619776
 5.984928986834,  5.139508800995,  1.871507160824
 6.070319282624,  3.906281101007,  3.443721502880
 5.146744668501,  2.543395351817,  4.579306095414
 3.722484798752,  2.925705980156,  5.930351026896
 2.517816229997,  4.503329314402,  5.685650690368
 2.827817886788,  5.800970314055,  4.195675251374
 4.290081136959,  5.706079435100,  2.834505934833
 4.168701850371,  4.449269813051,  1.283475102544
 2.906975871698,  2.900694574711,  1.383286354130
 3.009129841832,  1.576458192213,  2.878599889697
 1.808018019344, -0.022127483788,  2.921660618517
\end{verbatim}

\textbf{$7_6$: 30 vertices}

\begin{verbatim}
 0.324300143824, -0.318496069319,  4.509732543638
 0.683174040445,  0.832341161106,  6.105596947462
 1.986066871412,  2.292797026354,  5.693845885920
 2.999151071721,  3.017074759669,  4.128891519778
 4.736162553827,  4.006895140362,  4.073694373138
 6.054972518729,  3.668637181467,  2.608661325826
 5.674526913452,  2.302924026321,  1.197956282493
 4.331494063549,  0.930045879575,  1.756049785651
 3.563746598863, -0.311901619566,  3.122845395773
 2.230111005258,  0.377546611817,  4.444237555472
 0.706251527073,  1.658437632578,  4.251441072942
-0.393498666315,  3.290116153486,  3.893426106005
 0.371291667125,  4.782529577283,  2.803564684378
 2.268308349372,  4.654024489933,  2.183230782516
 2.612329462061,  2.901854764459,  1.282369679373
 3.490760823002,  1.606919859010,  0.036776674108
 5.159972010096,  0.508538354640, -0.048617024372
 6.287542413832,  0.409983919604,  1.600283217593
 5.699580564872,  1.594962385784,  3.100325516616
 4.218831206492,  2.807527578925,  2.519750809932
 4.043049259870,  4.341557144103,  1.248588349662
 2.383103722021,  4.576380548923,  0.157974403336
 0.796497161651,  3.787839703112,  1.085812189308
 1.316788028040,  2.986623448840,  2.842897824172
 1.351772641596,  4.154223734057,  4.466315715774
 2.829039876959,  4.135545630824,  5.814400521631
 4.118607574581,  2.607403540184,  5.772010616415
 4.056291313385,  1.289723457226,  4.268735564215
 2.696550342973,  1.497849593246,  2.816913259356
 1.080832896408,  0.319956726429,  2.771907402300
\end{verbatim}

\textbf{$7_7$: 30 vertices}

\begin{verbatim}
-0.076096572657,  1.401410357539,  1.322267310046
 0.157256423765,  3.014069693681,  2.481955304558
 1.153069105714,  4.508499382908,  3.362317165629
 1.649001649014,  5.973245215287,  4.630611515938
 3.346730535096,  7.003531343985,  4.393488926598
 4.768876793397,  6.063707346160,  3.347434017373
 4.271740741736,  4.224164356093,  2.740031898646
 5.579506430530,  2.807976198371,  3.273098503436
 5.672602324654,  2.628536574058,  5.262855919737
 4.431684410066,  3.915822607518,  6.158968299582
 3.447256009624,  5.019876613076,  4.812875629216
 2.818053805374,  5.615109669667,  3.010154853691
 2.452648163339,  6.850865676391,  1.480649767523
 1.972248975464,  6.007838611829, -0.268213288495
 1.652380373056,  4.040168651537, -0.107095305146
 1.622096517872,  2.503483678529,  1.172624266419
 1.317588494345,  1.378652469499,  2.798056398131
 2.738955936189,  1.472493297179,  4.201946883292
 3.902618950192,  3.089109839367,  4.382056142808
 5.361005644002,  4.456616300082,  4.437140906778
 6.213172347209,  4.812135720179,  2.663044558808
 4.923508987829,  5.489539500566,  1.292682366616
 3.010844819319,  4.935197944891,  1.107152805231
 1.058197141074,  5.193912047072,  1.453900635761
-0.633530331420,  4.883689182837,  2.474596444175
-0.580853234351,  3.842938707786,  4.181659293669
 1.065047798123,  2.713578649706,  4.306381587238
 2.482400780651,  3.035721231935,  2.932578588394
 3.160982515854,  1.431583944284,  1.949588053841
 1.671995535639,  0.480443480451,  1.012438337516
\end{verbatim}

\textbf{$8_1$: 32 vertices}

\begin{verbatim}
-2.008400903993, -0.277965842754,  0.525778844556
-0.401015277042,  0.572218140532,  1.358544502369
 0.866430418594,  2.088859263665,  1.664127735282
 2.837349695916,  1.832141689795,  1.441477135579
 3.542662930157,  0.109416144938,  0.710211152027
 2.175970943230, -0.860226107862, -0.381552061520
 0.389404667244,  0.037817260632, -0.422771593461
-1.133285499219,  1.332870274480, -0.357561561404
-2.158982786050,  1.556760252599,  1.344736401200
-1.967073566635, -0.045048105027,  2.526845041276
-1.320747343340, -1.753624949360,  1.712579986094
-2.001095122641, -2.528978210991, -0.000882653966
-0.696059372785, -2.988846098184, -1.444975177985
 0.555748216128, -1.460992784556, -1.759048437553
-0.288239485229,  0.245190518553, -2.372744857739
 0.932482841833,  1.449728695985, -3.401785566032
 2.466266504546,  2.024779133135, -2.254255753467
 2.307339853568,  1.219630106852, -0.430392445791
 1.589789612401,  0.272330260390,  1.178254393106
 0.805121528170,  0.573045090658,  2.993155654455
-0.563748215230,  2.026027891809,  3.115779490480
-0.834811805501,  3.286533099468,  1.586841574313
 0.378356748654,  3.127295266717,  0.004796893581
 0.645868387982,  1.762282830771, -1.432276024727
 1.697659173525,  0.174575524556, -2.042949139610
 0.918619366844, -0.778334037337, -3.619356557272
-0.887112560571, -1.507858357472, -3.164245062529
-1.362943744205, -1.096270812701, -1.265777159805
-0.110586358596, -1.841446175370,  0.104009194917
-0.609576328706, -3.523616530904,  1.063859133550
-2.483528654773, -3.398624313376,  1.751376133474
-3.281863864278, -1.629669119641,  1.268200786605
\end{verbatim}

\textbf{$8_{19}$: 27 vertices}

\begin{verbatim}
 0.000000000000,  0.000000000000,  0.000000000000
 2.000000000000,  0.000000000000,  0.000000000000
 3.732559997921,  0.999117537432,  0.000000000000
 3.740093269203,  2.983378029813,  0.250306908012
 2.836084633781,  4.664709871575, -0.346259622194
 1.080842092749,  5.230742730454, -1.120036350210
-0.586675178897,  4.705934378788, -0.148459105321
-1.433099231827,  3.361408147003,  1.066372378739
-0.986049984598,  1.802270274775,  2.236515987847
 0.971340008986,  1.400213471558,  2.153001078983
 2.033084618021,  1.944394769874,  0.547832832427
 2.320762514764,  2.900451360013, -1.185140688483
 2.528475469438,  4.521538634927, -2.337933547402
 2.738233262622,  6.323371275927, -1.495674203717
 1.935922742134,  6.339745551429,  0.336272769264
 0.946210849696,  5.175777804535,  1.626873193160
 0.338000172537,  3.315262219960,  2.037435650945
 0.017249180187,  2.072274549396,  0.503779740471
 1.000000000000,  1.374398953212, -1.092205998316
 2.838087418528,  1.050383933936, -1.810849796822
 4.241191111335,  2.464982055784, -1.637033281796
 4.848226791805,  4.221795228547, -0.898714137627
 4.369353897373,  4.765037706289,  0.965573490765
 2.607013977251,  4.027401374624,  1.557224273031
 1.140781868755,  3.728985312871,  0.230155601822
-0.055893490092,  3.090566067820, -1.239670562416
-1.020060848355,  1.346478391299, -1.070734237531
\end{verbatim}

\textbf{$10_{124}$: 32 vertices}

\begin{verbatim}
 0.000000000000,  0.000000000000,  0.000000000000
 2.000000000000,  0.000000000000,  0.000000000000
 3.749625987840,  0.968921515229,  0.000000000000
 4.260779271447,  2.111991727285,  1.559523264004
 4.723363038490,  3.242835425368,  3.142940072356
 5.074955844852,  2.393586002441,  4.919216491209
 3.823908832241,  0.939026480065,  5.484139466580
 2.016558949319,  0.111696348083,  5.262754612958
 0.994966724051, -0.147740463023,  3.563036440343
 0.530994529532,  0.698207939819,  1.811151426296
 0.945010187460,  1.768261727264,  0.172989063474
 2.216476928008,  1.482013102747, -1.344058744817
 3.260004049508, -0.211096385608, -1.554846332536
 3.714630934182, -1.032036987808,  0.211327759504
 4.009336867362,  0.131304000668,  1.811260932321
 4.748025325232,  1.237137538923,  3.305072014044
 3.239287762862,  2.280197843545,  4.102402666169
 1.576866177892,  2.074299324508,  5.195087558387
-0.088376878079,  1.351503600275,  4.355723569660
-0.935717586627,  0.095486260249,  3.050183771265
 0.181919811525, -1.262761590985,  2.098317488182
 2.061674393602, -0.603780542037,  1.918688180196
 2.416437938648,  1.312169603251,  1.467824710976
 2.921991742853,  2.790145448692,  0.218824009877
 4.848280879656,  2.629034119566, -0.294453011877
 5.776435193179,  1.138251111144,  0.662675815444
 5.958784709125, -0.009518512090,  2.290365449870
 4.575083550379, -0.755506009877,  3.526841796764
 2.933781670179,  0.386846731630,  3.492807225784
 1.512081032038,  1.764350925710,  3.207764262707
 0.024537393920,  2.611767941618,  2.173792151939
-0.962767246635,  1.569117884453,  0.781631814538
\end{verbatim}

\normalsize

Partial lattice configurations. First index is 1.

\tiny
\textbf{$6_1$: 26 vertices, 2-9, 16-23 on lattice}

\begin{verbatim}
 1.583389222588330,  0.295485149423139,  3.923108521431466
 2.828427124746190,  1.414213562373095,  2.828427124746190
 4.242640687119286,  2.828427124746190,  2.828427124746190
 4.242640687119286,  4.242640687119286,  1.414213562373095
 2.828427124746190,  4.242640687119286,  0.000000000000000
 1.414213562373095,  5.656854249492381,  0.000000000000000
 0.000000000000000,  5.656854249492381,  1.414213562373095
 0.000000000000000,  4.242640687119286,  2.828427124746190
 1.414213562373095,  2.828427124746190,  2.828427124746190
 1.000867245584298,  1.000867245584298,  2.129020440137794
 2.378269346607951, -0.444679645461330,  2.243727187471789
 3.530976312787990, -0.153795427197911,  3.852034689688109
 2.916925760140388,  1.481948501860211,  4.825319703105388
 1.016150987839091,  2.075182619639663,  4.637893593703873
-0.367101253559933,  1.996469075219529,  3.195528378306124
 0.000000000000000,  2.828427124746190,  1.414213562373095
 1.414213562373095,  4.242640687119286,  1.414213562373095
 2.828427124746190,  5.656854249492381,  1.414213562373095
 4.242640687119286,  5.656854249492381,  0.000000000000000
 5.656854249492381,  4.242640687119286,  0.000000000000000
 4.242640687119286,  2.828427124746190,  0.000000000000000
 2.828427124746190,  2.828427124746190,  1.414213562373095
 2.828427124746190,  4.242640687119286,  2.828427124746190
 2.491332739582196,  3.425590961867083,  4.622527769068077
 2.026453607801328,  2.471679819797600,  6.317797932160190
 1.481687497605007,  0.593616332331234,  5.898146383578224
\end{verbatim}

\textbf{$8_1$: 33 vertices, 1-3, 7-21, 26-33 on lattice}

\begin{verbatim}
 0.000000000000000,  1.414213562373095,  5.656854249492381
 0.000000000000000,  2.828427124746190,  4.242640687119286
 1.414213562373095,  4.242640687119286,  4.242640687119286
 2.828423589185952,  5.655340667830114,  4.308129346913210
 4.750944294310656,  6.153224620254956,  4.544826850669857
 6.126576831948576,  4.712356198457955,  4.722403066588766
 5.656854249492381,  2.828427124746190,  4.242640687119286
 4.242640687119286,  2.828427124746190,  2.828427124746190
 2.828427124746190,  4.242640687119286,  2.828427124746190
 1.414213562373095,  5.656854249492381,  2.828427124746190
 0.000000000000000,  5.656854249492381,  4.242640687119286
 0.000000000000000,  4.242640687119286,  5.656854249492381
 1.414213562373095,  2.828427124746190,  5.656854249492381
 1.414213562373095,  1.414213562373095,  4.242640687119286
 2.828427124746190,  1.414213562373095,  2.828427124746190
 2.828427124746190,  2.828427124746190,  1.414213562373095
 2.828427124746190,  4.242640687119286,  0.000000000000000
 4.242640687119286,  5.656854249492381,  0.000000000000000
 5.656854249492381,  5.656854249492381,  1.414213562373095
 5.656854249492381,  4.242640687119286,  2.828427124746190
 4.242640687119286,  4.242640687119286,  4.242640687119286
 3.798412823045159,  5.133244851387228,  5.977427655991470
 3.342546528170989,  7.051892552208125,  5.644295668136236
 3.482859476884064,  7.440403909802972,  3.687417746179106
 3.928917088119492,  5.882919598558162,  2.514696454616091
 4.242640687119286,  4.242640687119286,  1.414213562373095
 4.242640687119286,  2.828427124746190,  0.000000000000000
 2.828427124746190,  1.414213562373095,  0.000000000000000
 1.414213562373095,  1.414213562373095,  1.414213562373095
 1.414213562373095,  2.828427124746190,  2.828427124746190
 2.828427124746190,  2.828427124746190,  4.242640687119286
 2.828427124746190,  1.414213562373095,  5.656854249492381
 1.414213562373095,  0.000000000000000,  5.656854249492381
\end{verbatim}

\textbf{$10_{124}$: 32 vertices, 1-2, 6-11, 15-23, 28-32 on lattice}

\begin{verbatim}
 0.000000000000000,  0.000000000000000,  0.000000000000000
 1.414213562373095,  0.000000000000000, -1.414213562373095
 3.104706087393493, -0.276307247058916, -2.446632725344636
 4.632680037689242,  0.982214532378768, -2.161293387188179
 5.381990656518813,  2.613072669867294, -3.043809864031260
 4.242640687119286,  4.242640687119286, -2.828427124746190
 2.828427124746190,  4.242640687119286, -1.414213562373095
 2.828427124746190,  2.828427124746190,  0.000000000000000
 2.828427124746190,  1.414213562373095,  1.414213562373095
 1.414213562373095,  0.000000000000000,  1.414213562373095
 1.414213562373095, -1.414213562373095,  0.000000000000000
 2.137668947846121, -2.028073866201677, -1.760621433646020
 1.358050649012341, -0.886649269881901, -3.206078177737504
 1.760849488275634,  1.067605920886133, -3.069550728321591
 2.828427124746190,  1.414213562373095, -1.414213562373095
 4.242640687119286,  1.414213562373095,  0.000000000000000
 4.242640687119286,  2.828427124746190,  1.414213562373095
 2.828427124746190,  2.828427124746190,  2.828427124746190
 1.414213562373095,  4.242640687119286,  2.828427124746190
 0.000000000000000,  4.242640687119286,  1.414213562373095
 0.000000000000000,  2.828427124746190,  0.000000000000000
 1.414213562373095,  1.414213562373095,  0.000000000000000
 2.828427124746190,  0.000000000000000,  0.000000000000000
 4.439256267216421, -0.635214074108347, -1.000865902018410
 4.964895761766690, -0.838313653179917, -2.919837940050828
 4.028970632603468,  0.546900806951039, -4.017665363668194
 3.419943461040523,  2.236939072860518, -3.138542753299688
 4.242640687119286,  2.828427124746190, -1.414213562373095
 4.242640687119286,  4.242640687119286,  0.000000000000000
 2.828427124746190,  4.242640687119286,  1.414213562373095
 1.414213562373095,  2.828427124746190,  1.414213562373095
 0.000000000000000,  1.414213562373095,  1.414213562373095
\end{verbatim}
\normalsize
\section{Exact Edge Vector Construction}

The following procedure and table may be used to construct a non-interfering $5_1$ knot with equilateral edges. First define the $i$th edge vector using rational coefficients $a_i$ and $b_i$ from Table II:

\begin{equation}
    E_i(a_i,b_i)=\langle2\cos(\pi a_i)\cos(\pi b_i), 2\sin(\pi a_i)\cos(\pi b_i), 2\sin(\pi b_i)\rangle
\end{equation}

The first $N-2$ coordinates (19 in this case) can be constructed from the cumulative sum of these edges. We find the position of the 19th vertex:

\begin{equation}
R_{19}=-\sum_{i=1}^{19}E_i
\end{equation}

The last two edges if computed from the table will not close the knot properly. To ensure the knot is closed with an equilateral edge we define 

\begin{equation}
W=\frac{E_{20}-E_{21}}{2}
\end{equation}
and
\begin{equation}
    U=W-\frac{W\cdot R}{R\cdot R}R\qquad
\lambda=
\sqrt{
\frac{4-\frac14 R\cdot R}{U\cdot U}
},
\end{equation}

Then the modified final edges may be written:
\begin{equation}
E_{20}=\frac{R}{2}+\lambda U,
\qquad
E_{21}=\frac{R}{2}-\lambda U.
\end{equation}
Then each vertex may be found by placing the first on the origin and finding the cumulative sum for all others.

\begin{table}[ht]
\centering
\caption{Rational angle coefficients for the exact 21-vertex
$5_1$ construction. The final two sets may be used to satisfy the equilateral closure constraint.}
\label{tab:51-rational-angle-components}
\begin{tabular}{c r r r r}
\hline
$i$
& $a_i^{\mathrm{num}}$
& $a_i^{\mathrm{den}}$
& $b_i^{\mathrm{num}}$
& $b_i^{\mathrm{den}}$ \\
\hline
1  &  7789 & 9814 &   812 & 9033 \\
2  & -9548 & 9991 & -1590 & 9641 \\
3  &   911 & 1264 & -1047 & 6352 \\
4  &  9118 & 9385 &   413 & 4597 \\
5  & -4952 & 8253 &  1243 & 6499 \\
6  & -1069 & 5090 &   145 & 4586 \\
7  & -1207 & 3877 & -2699 & 8089 \\
8  & -2911 & 8182 &    91 & 2371 \\
9  &   128 & 3117 &  1307 & 7894 \\
10 &  3733 & 8221 &   233 & 3432 \\
11 &  3851 & 8964 & -1217 & 3984 \\
12 &   791 & 1411 &   197 & 6934 \\
13 &  2896 & 3281 &  1653 & 7487 \\
14 & -6956 & 8745 &    38 & 1345 \\
15 & -6094 & 9171 & -1097 & 3589 \\
16 & -4728 & 6863 &   667 & 9851 \\
17 & -1327 & 4808 &  1371 & 8285 \\
18 &  1161 & 9601 &   162 & 4207 \\
19 &   508 & 6619 & -2493 & 7474 \\
\hline
20 &  -65 & 2617 &   256 & 8083 \\
21 &   2886 & 7903 & 1889 & 9867 \\
\end{tabular}
\end{table}




\end{document}